\documentclass[12pt,a4paper,reqno]{amsart}
\usepackage[utf8]{inputenc}
\usepackage{amsmath}
\usepackage{amsfonts}
\usepackage{amssymb}
\usepackage{amsmath,bbm,amssymb,amsxtra}
\usepackage{enumerate}
\usepackage{caption}
\allowdisplaybreaks
\usepackage{epsfig}
\usepackage{ae}
    \def\Xint#1{\mathchoice
    {\XXint\displaystyle\textstyle{#1}}%
    {\XXint\textstyle\scriptstyle{#1}}%
    {\XXint\scriptstyle\scriptscriptstyle{#1}}%
    {\XXint\scriptscriptstyle\scriptscriptstyle{#1}}%
    \!\int}

    \def\XXint#1#2#3{{\setbox0=\hbox{$#1{#2#3}{\int}$ }
    \vcenter{\hbox{$#2#3$ }}\kern-.6\wd0}}
    \def\dashint{\Xint-}

\theoremstyle{definition}

\newtheorem{lemma}{Lemma}[section]

\newtheorem{proposition}[lemma]{Proposition}

\newtheorem{theorem}[lemma]{Theorem}

\newtheorem{corollary}[lemma]{Corollary}

\newtheorem{remark}[lemma]{Remark}

\newtheorem{definition}[lemma]{Definition}
\newtheorem{exam}[lemma]{Example}

\newcommand{\prop}[1]{\begin{proposition}\label{#1}
\sl }
\newcommand{\eprop}{\end{proposition}}
\newcommand{\thm}[1]{\begin{theorem}\label{#1}
\sl }
\newcommand{\ethm}{\end{theorem}}
\newcommand{\cor}[1]{\begin{corollary}\label{#1}
\sl }
\newcommand{\ecor}{\end{corollary}}

\newcommand{\lem}[1]{\begin{lemma}\label{#1}
\sl }
\newcommand{\elem}{\end{lemma}}
\newcommand{\defin}[1]{\begin{definition}\label{#1}
\sl }
\newcommand{\edefin}{\end{definition}}
\newcommand{\beqno}{\begin{eqnarray*}}
\newcommand{\eeqno}{\end{eqnarray*}}
\newcommand{\beqla}[1] {\begin {eqnarray}\label{#1}}
\def\eeq {\end {eqnarray}}
\newcommand{\beq}{\begin {eqnarray}}
\newcommand{\real}{{\mathbb R}}

\newcommand{\integer}{{\mathbb Z}}
\newcommand{\R}{{\mathbb R}}
\newcommand{\nanu}{{\mathbb N}}

\newcommand{\complex}{{\mathbb C}}

\newcommand{\res}{\mathcal R}
\newcommand{\nres}{{\mathcal R}^*}
\newcommand{\ext}{\mathcal E}
\newcommand{\next}{{\mathcal E}^*}
\newcommand{\integral}{{\mathbb I}}

\newcommand{\diam}{{\rm diam}\,}
\newcommand{\dist}{{\rm dist}\,}

\newcommand{\Lip}{{\rm Lip}\,}

\newcommand{\df}{\dot{\mathcal F}}
\newcommand{\ndf}{{\mathcal F}}
\newcommand{\dfp}{\dot{\mathcal F}^{s}_{p,q}(Z)}

\newcommand{\ndfp}{{\mathcal F}^{s}_{p,q}(Z)}
\newcommand{\db}{\dot{\mathcal B}}
\newcommand{\dbp}{\dot{\mathcal B}^{s}_{p,q}(Z)}
\newcommand{\ndb}{{\mathcal B}}
\newcommand{\ndbp}{{\mathcal B}^{s}_{p,q}(Z)}
\newcommand{\dm}{\dot{M}}
\newcommand{\dmp}{\dot{M}^s_{p,q}(Z)}

\newcommand{\dnp}{\dot{N}^s_{p,q}(Z)}

\newcommand{\J}{{\mathcal J}}
\newcommand{\I}{{\mathcal I}}

\newcommand{\jpe}{\J^{s}_{p,q}(E)}

\newcommand{\ipe}{\I^{s}_{p,q}(E)}

\newcommand{\locint}{L^1_{\rm{loc}}}

\title[Traces of function spaces]{Traces of Besov, Triebel-Lizorkin and Sobolev spaces on metric spaces}

\author[Saksman]{Eero Saksman}
\address{Department of Mathematics and Statistics, University of Helsinki, PO~Box~68, FI-00014 Helsinki, Finland.}
\email{eero.saksman@helsinki.fi}
\thanks{E.S.\ was supported by the Finnish CoE in Analysis and Dynamics Research, and by the Academy of Finland, projects 113826 and 118765.}

\author[Soto]{Tom\'as Soto}
\address{Department of Mathematics and Statistics, University of Helsinki, PO~Box~68, FI-00014 Helsinki, Finland}
\email{tomas.soto@helsinki.fi}
\thanks{T.S.\ was supported by the Finnish CoE in Analysis and Dynamics Research, and by the V\"ais\"al\"a foundation.}



\begin{document}

\maketitle

\begin{abstract}
We establish trace theorems for function spaces defined on general Ahlfors regular metric spaces $Z$. The results cover the Triebel-Lizorkin spaces and the Besov spaces for smoothness indices $s<1,$ as well as the first order Haj\l asz-Sobolev space $M^{1,p}(Z)$. They generalize the classical results from the Euclidean setting, since the traces of these function spaces onto any closed Ahlfors regular subset $F \subset Z$ are Besov spaces defined intrinsically on $F$. Our method employs the definitions of the function spaces via hyperbolic fillings of the underlying metric space.
\end{abstract}

\section{Introduction}\label{se:introduction}

A classical fact, originally due to Gagliardo \cite{Ga}, states that the traces of the Sobolev space $W^{1,p}(\R^d)$, $p\in (1,\infty)$, on the hyperplane $\R^{d-1}\times\{0\}$ lie in the Besov space $B^{1-1/p}_{p,p}(\R^{d-1})$ and, conversely, any function in $B^{1-1/p}_{p,p}(\real^{d-1})$ is a trace of some function in $W^{1,p}(\R^d)$. This important result has been generalized to many other function spaces, most notably to the Triebel-Lizorkin spaces ${F}^s_{p,q}$ and the Besov spaces ${B}^s_{p,q}$. Loosely speaking, we have
\[
  {B}^s_{p,q}(\real^d)_{|\real^{d-1}} = {B}^{s-1/p}_{p,q}(\real^{d-1}) \quad \text{and} \quad {F}^s_{p,q}(\real^d)_{|\real^{d-1}} = {B}^{s-1/p}_{p,p}(\real^{d-1})
\]
for $p \geq 1$ and $s > 1/p$. We refer \cite{Ni,J2,Pe,T} and the references therein for these facts, as well as generalizations to some classes of subdomains of $\real^d$.

The Besov spaces, and later on the Triebel-Lizorkin spaces, have been studied in the fairly general setting of doubling metric measure spaces; we refer to e.g.~\cite{HS,BP,HMY,GKZ,KYZ} and the references therein, although this list is by no means comprehensive. Especially the full scales of these spaces in the setting of doubling metric spaces were introduced in \cite{KYZ}, and in this paper we shall work with the equivalent definitions given in \cite{BSS,S} in terms of the "hyperbolic fillings" of the underlying metric space -- the actual definitions are given in the next section.

In order to describe our results, let $Z := (Z,d,\mu)$ be a $Q$-Ahlfors regular metric measure space for some $Q > 0$, and let $F \subset Z$ be a closed $\lambda$-Ahlfors regular subset, where $\lambda \in (0,Q]$. We equip $F$ with the metric $d_{|F}$ and the Hausdorff $\lambda$-measure.

\thm{th:besov-trace} Let $F\subset Z$ be a closed $\lambda$-Ahlfors regular subset.
Suppose that $0 < s < 1$, $\max\big(Q/(\lambda+s) , (Q-\lambda)/s\big) < p < \infty$ and $0 < q \leq \infty$. Then there exist bounded linear operators
\[
  \res\colon\dbp \to \db^{s-\frac{Q-\lambda}{p}}_{p,q}(F) \quad \text{and} \quad \ext\colon \db^{s-\frac{Q-\lambda}{p}}_{p,q}(F)\to\dbp
\]
such that

\smallskip
(i) $\res f = f_{|F}$ for all continuous functions $f$ in $\dbp$, and

(ii) $\res \big(\ext f) = f$ for all $f \in \db^{s-\frac{Q-\lambda}{p}}_{p,q}(F)$.
\ethm

We refer to Remark \ref{re:assumptions} below for a concrete explanation of the range of the parameter $p$, as well as an alternative way to interpret part (i) of the statement. A similar result in the range $p > 1$  and $q\geq 1$ has very recently been obtained in \cite{Mar} using interpolation techniques.

To formulate the two other trace theorems, we need a minor additional condition on the subset $F$. 

\begin{definition} The closed set $F\subset Z$ is \emph{porous} if there exists a constant $c \in (0,1)$ such that for all balls $B \subset Z$ with radius $r < \diam Z$ such that $B \cap F \neq \emptyset$, there exists $\xi \in Z$ such that $B(\xi,cr) \subset B \backslash F$. 
\end{definition}

\begin{remark}\label{re:porous} For porous $\lambda$-Ahlfors regular sets $F\subset Z$, it follows that $\lambda<Q$; see e.g.~\cite[Proposition 3.4]{JJKRRS}. On the other hand, if $F \subset Z$ is $\lambda$-Ahlfors regular with $\lambda < Q$, it follows that $F$ is porous subset of $Z$ \cite[Theorem 5.3]{JJKRRS}.
\end{remark}

Our trace theorem for the Triebel-Lizorkin spaces reads as follows.

\thm{th:triebel-trace} Let $F\subset Z$ be a closed and porous $\lambda$-Ahlfors regular subset.
Suppose that $0 < s < 1$, $\max\big(Q/(\lambda+s) , (Q-\lambda)/s\big) < p < \infty$ and $Q/(Q+s) < q \leq \infty$. Then there exist bounded linear operators
\[
  \res\colon\dfp \to \db^{s-\frac{Q-\lambda}{p}}_{p,p}(F) \quad \text{and} \quad \ext\colon \db^{s-\frac{Q-\lambda}{p}}_{p,p}(F)\to\dfp
\]
such that

\smallskip
(i) $\res f = f_{|F}$ for all continuous functions $f$ in $\dfp$, and

(ii) $\res \big(\ext f) = f$ for all $\db^{s-\frac{Q-\lambda}{p}}_{p,p}(F)$.
\ethm

Finally, our trace theorem for the Haj\l asz-Sobolev spaces $\dm^{1,p}$ reads as follows. A similar but much weaker result was established in \cite{GKS}.

\thm{th:sobolev-trace} Let $F\subset Z$ be a closed and porous $\lambda$-Ahlfors regular subset.
Suppose that $\max\big(Q/(\lambda+1) , Q-\lambda\big) < p < \infty$. Then there exist bounded linear operators
\[
  \res\colon\dm^{1,p}(Z) \to \db^{1-\frac{Q-\lambda}{p}}_{p,p}(F) \quad \text{and} \quad \ext\colon \db^{1-\frac{Q-\lambda}{p}}_{p,p}(F)\to\dm^{1,p}(Z)
\]
such that

\smallskip
(i) $\res f = f_{|F}$ for all continuous functions $f$ in $\dm^{1,p}(Z)$, and

(ii) $\res \big(\ext f) = f$ for all $\db^{1-\frac{Q-\lambda}{p}}_{p,p}(F)$.
\ethm

One should observe that Theorems \ref{th:besov-trace}, \ref{th:triebel-trace} and \ref{th:sobolev-trace} are exactly of the same form as the classical results in the Euclidean setting, and that Theorems \ref{th:triebel-trace} and \ref{th:sobolev-trace} (along with their non-homogeneous counterparts below) are completely new in this generality. A (very incomplete) list of previous results in the setting where $F$ is a subset of an Euclidean space includes \cite{JoW, Sa1, Sa2, T4, IV, Jo, CaHa,HaMa} -- the trace spaces appearing in these papers are sometimes defined in a non-intrinsic manner, however e.g. \cite{Sa1, Sa2,IV} employ intrinsic characterizations in terms of optimal polynomial approximations. Recently a trace theorem similar to Theorem \ref{th:besov-trace} for BV functions in the setting of metric measure spaces was also obtained in \cite{MSS}. We further refer to \cite{HIT} for metric results concerning the restriction and extension of Besov and Triebel-Lizorkin functions to subsets which are
sufficiently ``thick'', i.e.~have positive $\mu$-measure and satisfy a certain measure density condition.

We point out that something like porosity needs to be assumed in Theorems \ref{th:triebel-trace} and \ref{th:sobolev-trace} -- $F$ can not be properly $Q$-dimensional e.g.~in terms of the measure density condition considered in \cite{HIT}.

\begin{exam}\label{ex:1} Let us give a simple application of our results illustrating a curious phenomenon concerning smoothness spaces defined on fractal subsets of an Euclidean space.

Consider a self-similar fractal subset $Z$ of $\real^d$ in the sense of Hutchinson \cite{Hu} generated by a collection of similitudes $S_i\colon \real^d\to\real^d$, $1 \leq i \leq N$, satisfying the so-called open set condition. Let $F$ be a sub-fractal generated by a proper subcollection of $(S_i)_{1 \leq i \leq N}$. Then by \cite[Theorem 4.14]{Mat} and \cite[Theorem 5.3]{JJKRRS}, $Z$ and $F$ are respectively $Q$-Ahlfors regular and $\lambda$-Ahlfors regular with $0 < \lambda < Q < d$, $Z$ is a porous subset of $\real^d$ and $F$ is a porous subset of both $Z$ and $\real^d$.

In particular, our results imply that the function space $\db^{\sigma}_{p,p}(F)$ with $\lambda \max(0,p^{-1}-1) < \sigma \leq 1 - (d-\lambda)/p$ can be realized as the trace space of an appropriate Sobolev or Triebel-Lizorkin space defined either on $\real^d$ or on $Z$!
\end{exam}
It is further natural to ask whether analogous trace results hold for non-homogeneous versions of these function spaces. In Section \ref{se:nonhomogeneous} we  define the non-homogeneous function spaces $\ndbp$ and $\ndfp$, establish some of their basic properties and show that we have the following counterpart for our homogeneous trace theorems.

\thm{th:nonhomogeneous}
The Theorems \ref{th:besov-trace}, \ref{th:triebel-trace} and \ref{th:sobolev-trace} hold with $\ndb$, $\ndf$ and $M$ in place of $\db$, $\df$ and $\dm$ respectively.
\ethm

We finally point out that the Ahlfors regularity of the spaces $Z$ and $F$ is not strictly needed in these results, and a property known in previous literature as \emph{Ahlfors co-regularity} of $F$ with respect to $Z$ would suffice. Section \ref{se:appendix} is an appendix where we elaborate on this, as well as on some technicalities that are needed in the proofs of our main results.

\section{Preliminaries}\label{se:preliminaries}

In this section we give the definitions of the relevant function spaces and some details concerning them. In the rest of this section, $Z := (Z,d,\mu)$ is assumed to be a doubling metric measure space such that the measure $\mu$ is Borel regular and every ball $B(\xi,r) := \{\eta \in Z \,:\, d(\eta,\xi) < r\}$ has positive and finite $\mu$-measure. The \emph{doubling} assumption means that there exists a constant $c \in [1,\infty)$ such that $\mu\big(B(\xi,2r)\big) \leq c \mu\big(B(\xi,r)\big)$ for all $\xi \in Z$ and $r > 0$. It follows from this assumption that there exist constants $C \geq 1$ and $Q > 0$ such that the measure $\mu$ satisfies
\beqla{eq:doubling}
  \mu\big(B(\xi,\lambda r) \big) \leq C \lambda^Q \mu\big(B(\xi,r) \big)
\eeq
for all $\xi \in Z$, $r > 0 $ and $\lambda \geq 1$. $Q$ is in a sense an upper bound for the dimension of $Z$, and it will be fixed from now on.

Let us introduce some notation conventions that will be used in this section, as well as in the later sections with obvious modifications. For an arbitrary ball $B \subset Z$ with a distinguished center point $\xi \in Z$ and radius $r > 0$, we write $\lambda B := B(\xi,\lambda r)$ for $\lambda > 0$. If $f$ is a complex-valued function on $Z$ and $E \subset Z$, we write
\[
  \dashint_E f d\mu := \frac{1}{\mu(E)} \int_E f d\mu
\]
whenever the latter quantity is well-defined. The notation $\locint(Z)$ will (instead of the usual one) stand for the space of complex-valued $\mu$-measurable functions on $Z$ that are integrable on every ball $B \subset Z$. Finally, we will use the notations $\lesssim$, $\gtrsim$ and $\approx$ when dealing with unimportant multiplicative constants. More precisely, when $f$ and $g$ are non-negative functions with the same domain, the notations $f \lesssim g$ or $g \gtrsim f$ mean that there exists a positive constant $c$, usually independent of some parameters obvious from the context, such that $f \leq c g$ on the domain of $f$ and $g$. The notation $f \approx g$ means that $f \lesssim g$ and $g \lesssim f$.

The construction referred to above as the \emph{hyperbolic filling} of $Z$ is roughly speaking a graph $(X,E)$ such that if $Z$ is nice enough and $(X,E)$ is endowed with its natural path metric, $(X,E)$ is hyperbolic in the sense of Gromov and its boundary at infinity coincides with $Z$. We refer to the introduction section of \cite{BSS} for a detailed explanation of the motivation of this construction in the context of function spaces. We shall next explain the actual construction of $(X,E)$.

For all $n \in \integer$, let $(\xi_x)_{x \in X_n}$, where $X_n$ is a suitable index set, be a maximal set of points in $Z$ such that $d(\xi_x,\xi_{x'}) \geq 2^{-n-1}$ for all pairwise distinct $x$, $x' \in X_n$. Write $B(x) := B(\xi_x,2^{-n})$ for all $x \in X_n$. It is easily seen that the balls $2^{-1} B(x)$, $x \in X_n$, cover $Z$, and the doubling assumption implies that the balls $B(x)$, $x \in X_n$, have bounded overlap (uniformly in $n$). Write $|x| := n$ for all $x \in X_n$ (the ``level'' of $x$). We then consider the disjoint union $X := \bigsqcup_{n \in \integer} X_n$, and denote by $(X,E)$ the graph such that the vertices $x$, $x' \in X$ are joined by an edge in $E$ if and only if $x \neq x'$, $||x|-|x'|| \leq 1$ and $B(x)\cap B(x') \neq \emptyset$; in this case we write $x \sim x'$. The ``Poisson extension'' $Pf\colon X\to \complex$ of a function $f \in \locint(Z)$ is defined by
\[
  Pf(x) := \dashint_{B(x)} f d\mu
\]
for all $x \in X$.

We equip the edges in $E$ with an orientation and denote by $e_{x,x'}$ the directed edge from $x$ to $x'$ for any two neighbors $x$ and $x'$ in $X$. The orientation is chosen so that if $x \sim x'$ and $|x| < |x'|$, then $x'$ is the endpoint of the edge joining $x$ and $x'$. For an edge $e \in E$, denote by $e_-$ the starting point and by $e_+$ the endpoint of $e$. For a sequence $u\colon X \to \complex$, we define the discrete derivative $du \colon E \to \complex$ by $du(e) = u(e_+) - u(e_-)$ for all $e \in E$. Finally, we write $|e| := \min\big(|e_-|,|e_+|\big)$ (the ``level'' of $e$) and $B(e) := B(e_-)\cup B(e_+)$ for all $e \in E$.

The definitions of the spaces $\db^s_{p,q}$, $\df^s_{p,q}$ and $\dm^{1,p}$ then read as follows. 

\defin{de:spaces}
{\rm (i)} Let $0 < s \leq 1$, $Q/(Q+s) < p \leq \infty$ and $0 < q \leq \infty$. Then $\ipe$ is the quasi-normed space of sequences $u\colon E\to \complex$ such that
\beqla{eq:besov-norm}
  \| u \|_{\ipe} := \bigg( \sum_{k \in \integer} 2^{ksq} \big\| \sum_{|e| = k} |u(e)| \chi_{B(e)} \big\|_{L^p(Z)}^q \bigg)^{1/q}
\eeq
(standard modification for $q = \infty$) is finite. Furthermore, the homogeneous Besov space $\dbp$ is the quasi-normed space of functions $f \in \locint(Z)$ such that
\[
  \|f\|_{\dbp} := \big\| d(Pf) \big\|_{\ipe}
\]
is finite.

\smallskip
{\rm (ii)} Let $0 < s \leq 1$, $Q/(Q+s) < p < \infty$ and $Q/(Q+s) < q \leq \infty$. Then $\jpe$ is the quasi-normed space of sequences $u\colon E\to \complex$ such that
\beqla{eq:triebel-norm}
  \| u \|_{\jpe} := \bigg( \int_{Z} \Big( \sum_{e \in E} \big[2^{|e|s} |u(e)| \big]^q \chi_{B(e)}(\xi)  \Big)^{p/q} d\mu(\xi) \bigg)^{1/p}
\eeq
(standard modification for $q = \infty$) is finite. Furthermore, the homogeneous Triebel-Lizorkin space $\dfp$ is the quasi-normed space of functions $f \in \locint(Z)$ such that
\[
  \|f\|_{\dfp} := \big\| d(Pf) \big\|_{\jpe}
\]
is finite.

\smallskip
{\rm (iii)} Let $0 < s \leq 1$ and $0 < p < \infty$. The homogeneous Haj\l asz-Sobolev space $\dm^{s,p}(Z)$ is defined as the class of $\mu$-measurable functions $f\colon Z \to \complex$ such that there exists a function $g\colon Z \to [0,\infty]$ in $L^p(Z)$ such that
\[
  |f(\xi) - f(\eta)| \leq d(\xi,\eta)^s \big( g(\xi) + g(\eta) \big)
\]
for all $\xi$, $\eta \in Z$. The quasi-norm $\|f\|_{\dm^{s,p}(Z)}$ of a function $f \in \dm^{s,p}$ is obtained as the infimum of $\|g\|_{L^p}$ over all admissible $g$.
\edefin

\begin{remark}\label{re:space-properties}
{\rm (i)} Strictly speaking the spaces $\dbp$, $\dfp$ and $\dm^{s,p}(Z)$ become quasi-normed spaces after dividing out the functions $f$ such that $\|f\| = 0$, i.e.~the functions that are constant $\mu$-almost everywhere. In the sequel we shall abuse notation by writing $f \in \dbp$ for both functions $f$ and equivalence classes $f$ satisfying $\|f\|_{\dbp} < \infty$, and similarly for the other two families of function spaces introduced above. The precise meaning will be obvious from context.

\smallskip
{\rm (ii)} The spaces $\dfp$ and $\dbp$ were introduced in \cite{BSS} and \cite{S} respectively. We refer to these papers for all basic properties concerning these spaces. Let us only mention here that these two function spaces are quasi-Banach spaces for all admissible values of the parameters, and reflexive Banach spaces for $1 < p,\,q < \infty$. While the sequence spaces $\jpe$ and $\ipe$ obviously depend on the choice of the hyperbolic filling $(X,E)$, the spaces $\dfp$ and $\dbp$ do not -- we refer to Remark \ref{re:filling-choice} below for more information.

\smallskip
{\rm (iii)} The Besov spaces $\dnp$ and the Triebel-Lizorkin spaces $\dmp$ were introduced in \cite{KYZ} in the generality of all metric measure spaces. Under our assumptions and in the parameter ranges given in the definition above, they coincide with $\dbp$ and $\dfp$ respectively; see \cite[Propostion 3.1]{BSS} and \cite[Proposition 3.1]{S}. In particular, $\db^s_{p,q}(\real^d)$ and $\df^s_{p,q}(\real^d)$ with $0 < s < 1$ coincide with the standard Fourier-analytically defined Besov and Triebel-Lizorkin spaces on $\real^d$ for all admissible values of the parameters. The spaces $\dnp$ and $\dmp$ with $s \geq 1$ are often trivial \cite[Theorem 4.1]{GKZ}, but there is one exception important for us: $\df^1_{p,\infty}(Z) = \dm^1_{p,\infty}(Z) = \dm^{1,p}(Z)$, where $\dm^{1,p}$ is the standard first order Haj\l asz-Sobolev space, for $Q/(Q+1) < p < \infty$.

\smallskip
{\rm (iv)} The spaces $\dm^{1,p}(Z)$ were introduced in \cite{H}; see also \cite{H2}. For $1 < p < \infty$ they are one of the more well-known generalizations of the standard Sobolev spaces to the setting of metric measure spaces. In \cite{KS} it was shown that $\dm^{1,p}(\real^d)$ for $d/(d+1) < p \leq 1$ coincides with the homogeneous \emph{Hardy-Sobolev} space with the same indices.

\smallskip
{\rm (v)} We have the restriction $p > Q/(Q+s)$ in the definitions of the spaces $\dbp$ and $\dfp$. This is because the definitions of these spaces require a priori local integrability, and generally speaking it is for $p > d/(d+s)$ that the Besov and Triebel-Lizorkin distributions on $\real^d$ are locally integrable functions. We also impose a similar restriction on the parameter $q$ for the spaces $\dfp$ because of certain technical reasons which are common in the study of these spaces; see e.g.~the proofs in \cite{BSS} or \cite{T} for more information.

\smallskip
{\rm (vi)} It will be useful to note that we have the following equivalent quasinorm on the space $\ipe$ for all admissible parameters:
\beqla{eq:besov-norm-2}
  \| u \|_{\ipe} \approx \bigg( \sum_{k \in \integer} 2^{ksq} \Big( \sum_{|e| = k} \mu\big( B(e)\big)|u(e)|^p  \Big)^{q/p} \bigg)^{1/q}
\eeq
(obvious modifications for $p = \infty$ and/or $q = \infty$). This follows easily from \eqref{eq:besov-norm} and the fact that the sets $B(e)$, $|e| = k$, have bounded overlap uniformly in $k \in \integer$.

\smallskip
{\rm (vii)} Let $(A_e)_{e \in E}$ be a collection of measurable subsets of $Z$ such that $A_e \subset \big( \lambda B(e_-)\big) \cup \big( \lambda B(e_+)\big)$ for some uniform $\lambda \geq 1$ and $\inf_{e \in E} \mu(A_e)/\mu\big(B(e)\big) > 0$. Then we get equivalent quasinorms on $\ipe$ and $\jpe$ by replacing $\chi_{B(e)}$ by $\chi_{A_e}$ in \eqref{eq:besov-norm} and \eqref{eq:triebel-norm} respectively. This can be proven by a standard maximal function argument; see \cite[Proposition 2.2]{BSS}
\end{remark}

\begin{remark}\label{re:filling-choice}
As mentioned above, the choice of the hyperbolic filling $(X,E)$ is not unique, but this has no essential bearing on the classes $\dbp$ and $\dfp$ or their quasi-norms -- any two admissible choices yield equivalent quasi-norms for both spaces, with the equivalence constants independent of these two choices. However, in this paper we shall need even more flexibility in the choice of $(X,E)$ -- we want to choose the hyperbolic filling of $Z$ in such a way that a hyperbolic filling of a fixed subspace $F$ is obtained in a natural way as the ``restriction'' of $(X,E)$ to the edges corresponding to balls that lie ``above'' $F$. This choice is formulated as Lemma \ref{le:filling-choice} below.

To elaborate on the admissible flexibility, it is enough that we have $d(x,x') \geq c_1 2^{-n}$ for all distinct $x$, $x' \in X_n$, that the radii $r_x$ of the balls $B(x) := B(\xi_x,r_x)$ ($x \in X_n$) are comparable to $2^{-n}$ uniformly in $n$, and that the balls $\big( c_2 B(x) \big)_{x \in X_n}$ cover $Z$; here the constants $c_1 > 0$ and $c_2 \in (0,1)$ are uniform in $n$. Then $(X,E)$ can be constructed exactly as explained above, and the resulting spaces $\dbp$ and $\dfp$ have quasinorms essentially independent of the precise choice of the hyperbolic filling. We refer to \cite[Remark 2.8]{BSS} and \cite[Remark 2.8]{S} for details.
\end{remark}

\lem{le:filling-choice}
Suppose that $F$ is a closed subset of $Z$ equipped with the metric $d_{|F}$. Let $(X^Z,E^Z)$ be an admissible hyperbolic filling of $Z$, and write $X^F$ for the set of vertices $x \in X^Z$ such that $B(x) \cap F \neq \emptyset$. Then $(X^Z,E^Z)$ can be chosen such that the following properties hold.

\smallskip
{\rm (i)} The balls corresponding to the vertices in $X^F$ are centered in $F$, i.e.~$\xi_x \in F$ for all $x \in X^F$;

\smallskip
{\rm (ii)} The balls $\big(B(x)_{|F}\big)_{x \in X^F}$, generate an admissible hyperbolic filling $(X^F,E^F)$ of the metric space $(F,d_{|F})$.

\smallskip
In particular, we have $X^F \subset X^Z$ and $E^F \subset E^Z$ in a natural way.
\elem

\begin{proof}
For $n \in \integer$, let $(\xi_x)_{x \in X'_n}$ be a maximal $2^{-n-1}$-separated subset of $\{\xi \in Z \,:\, \dist(\xi,F) \geq 2^{-n}\}$, where $X'_n$ is a suitable index set. Furthermore, let $(\xi_x)_{x \in X''_n}$ be a maximal $2^{-n}$-separated subset of $F$, with again $X''_n$ a suitable index set. Let $B(x) := B(\xi_x,2^{-n}) \subset Z$ for $x \in X'_n$ and $B(x) := B(\xi,2^{-n+2}) \subset Z$ for $x \in X''_n$. Writing $X^Z_n := X'_n \cup X''_n$ and $X^Z := \bigsqcup_{n \in \integer} X^Z_n$, we are in the situation of Remark \ref{re:filling-choice} (with $c_1 = c_2 = 2^{-1}$), so the resulting graph $(X^Z,E^Z)$ is an admissible hyperbolic filling of $Z$.

In addition, for $x \in X^Z_n$ we have $B(x) \cap F \neq \emptyset$ if and only if $x \in X''_n$ (hence also $\xi_x \in F)$, and the balls $B(x)_{|F} \subset F$ corresponding to the vertices $x \in X''_n$ obviously generate an admissible hyperbolic filling $(X^F,E^F)$ for the metric space $(F,d_{|F})$.
\end{proof}
 
\section{Traces of Besov spaces}\label{se:besov}

In this section, we shall work with the assumptions of Theorem \ref{th:besov-trace}. In other words, the metric measure space $(Z,d,\mu)$ is assumed to be $Q$-Ahlfors regular (with $Q$ as in \eqref{eq:doubling}), which means that the measure $\mu$ satisfies $\mu\big(B(\xi,r)\big) \approx r^Q$ uniformly in $\xi \in Z$ and $0 < r < \diam Z$. We assume that $F$ is a closed subset of $Z$ of Hausdorff dimension $\lambda \in (0,Q]$, equipped with the metric $d_{|F}$. We denote its $\lambda$-Hausdorff measure by $\nu$,  and assume  it to be $\lambda$-Ahlfors regular. Write $(X^Z,E^Z)$ for the hyperbolic filling of $Z$, and similarly with $F$ in place of $Z$. The hyperbolic fillings are chosen so that $X^F_n$ is in a natural way a subset of $X^Z_n$ for all $n \in \integer$; see Lemma \ref{le:filling-choice} above.

To make sense of the trace spaces of the Besov spaces $\dbp$, let us begin by recalling some very basic properties of locally integrable functions.

For a measurable function $f\colon Z \to \complex$, denote by $\Lambda_f$ the set of points $\xi \in Z$ such that there exists a number $c_{\xi,f} \in \complex$ so that
\[
  \lim_{r \to 0} \dashint_{B(\xi,r)} \big|f - c_{\xi,f} \big| d\mu = 0.
\]
It is well known (\cite[Theorem 2.7]{He}) that the doubling property \eqref{eq:doubling} implies that if $f \in \locint(Z)$, then $\Lambda_f$ has full $\mu$-measure (namely it contains the Lebesgue points of $f$) and that it does not depend on the precise representative of $f$ (with respect to equality $\mu$-almost everywhere). The point is that $f$ as a function is essentially well-defined in $\Lambda_f$, and that under a fractional smoothness assumption on $f$, the set $Z\setminus \Lambda_f$ turns out to have a relatively small Hausdorff dimension. This is quantified in the following lemma, which is well known in the Euclidean setting, so we have only included an outline of a proof that is easily adapted to our setting.

\lem{le:hausdorff}
Suppose that $f$ is a function in $\dfp$ with $0 < s \leq 1$, $Q/(Q+s) < p \leq Q/s$ and $Q/(Q+s) < q \leq \infty$, or $f\in\dbp$ with $0 < s < 1$, $Q/(Q+s) < p \leq Q/s$ and $0 < q \leq \infty$. Then the Hausdorff dimension of $\Lambda\backslash f$ is at most $Q - ps$.
\elem

\begin{proof}
We first consider the case of Triebel-Lizorkin functions. With the parameters as in the statement, we have $\dfp \subset \df^s_{p,\infty}(Z) = \dm^{s,p}(Z)$ \cite[Proposition 3.1]{BSS}, so it suffices to verify the statement for the latter space. Let $g \in L^p(Z)$ be a Haj\l asz $s$-gradient of a function $f \in \dm^{s,p}(Z)$, and fix $\epsilon \in(0,s)$. Taking $\xi \in Z$ and $0 < r_1 < r_2 < 1$, and $k \in \nanu_0$ such that $2^k r_1 < r_2 \leq 2^{k+1}r_1$, the doubling condition and the weak $(1,p)$-Poincar\'e inequality satified by the functions of $\dm^{s,p}(Z)$ (see \cite[Theorem 8.7]{H2} and \cite[Lemma 4.1]{KYZ2}) yield
\begin{align*}
\big|f_{B(\xi,r_2)} - f_{B(\xi,r_1)}\big|
& \leq \sum_{n=0}^{k} \big| f_{B(\xi,2^{n+1}r_1)} - f_{B(\xi,2^nr_1)} \big| + \big|f_{B(\xi,2^{k+1}r_1)} - f_{B(\xi,r_2)} \big| \\
& \lesssim \sum_{n=0}^{k+1} \big(2^n r_1\big)^{\epsilon} \sup_{r \in (0,4)} \bigg(r^{(s-\epsilon)p-Q} \int_{B(\xi,r)} g^p d\mu \bigg)^{1/p} \\
& \approx (r_2)^\epsilon \sup_{r \in (0,4)} \bigg(r^{(s-\epsilon)p-Q} \int_{B(\xi,r)} g^p d\mu \bigg)^{1/p}.
\end{align*}
Now $\xi \in Z\backslash \Lambda_f$ only if the latter supremum is infinite (otherwise one can take $c_{\xi,f} = \lim_{r\to 0} \dashint_{B(\xi,r)} f d\mu$), and since $g \in L^p$, a standard covering argument shows that this happens in a set of Hausdorff $(Q-sp+\epsilon p)$-content zero, so the Hausdorff dimension of $Z\backslash \Lambda_f$ is at most $Q-sp + \epsilon p$. Letting $\epsilon \to 0$ yields the desired upper bound.

Let us now consider the case of Besov spaces. Suppose that the parameters are as in the statement, and take $\epsilon \in (0,s)$ arbitrarily close to $s$. We have $\dbp = \dnp$  (see \cite{KYZ} for the definition of the latter space), and it is easily seen that for any $f \in \dnp$, we have $f \in \dm^{\epsilon,p}(B)$ for all balls $B \subset Z$. By the first part of the proof, this means that the Hausdorff dimension of $Z\backslash \Lambda_f$ is at most $Q - \epsilon p$, and taking $\epsilon \to s$ yields the desired upper bound.
\end{proof}

\begin{remark}\label{re:hausdorff}
When $p > Q/s$, the functions $f$ in $\dfp$ and $\dbp$ coincide with  H\"older continuous functions $\mu$-almost everywhere, which means that the set $Z\backslash\Lambda_f$ is empty. For the spaces $\dfp$, this follows again from the embedding $\dfp \subset \dm^{s,p}$ and a Haj\l asz's Sobolev-type embedding theorem for the spaces $\dm^{s,p}$ (\cite[Theorem 8.7]{H2} and \cite[Lemma 4.1]{KYZ2}). For the spaces $\dbp$, we may again note that the functions in $\dnp$ are locally in $\dm^{\epsilon,p}$ for all $\epsilon \in (0,s)$ such that $p > Q/\epsilon$.

We also note that the Ahlfors regularity of $Z$ is not strictly speaking needed here; a closer examination of the proof shows that the doubling condition \eqref{eq:doubling} suffices.
\end{remark}

With this in mind, we are in a position to give the proof of Theorem \ref{th:besov-trace}. Before the proof, we still make some remarks concerning the actual statement and define a collection of auxiliary functions.

\begin{remark}\label{re:assumptions}
Let us elaborate on the precise assumptions and the statement of Theorem \ref{th:besov-trace}. Firstly, the condition $p > Q/(\lambda+s)$ is equivalent with the requirements that $p > Q/(Q+s)$ and $p > \lambda/(\lambda+(s-(Q-\lambda)/p)$, so that the functions in the Besov spaces in question can be expected to be locally integrable in the first place. The condition $p > (Q-\lambda)/s$ comes from the requirement that $s - (Q-\lambda)/p > 0$, and by Lemma \ref{le:hausdorff} and Remark \ref{re:hausdorff}, this also means that the functions in $\dbp$ are essentially well-defined $\nu$-almost everywhere in $F$, and in this way part (i) of the statement also makes sense for all (not necessarily continuous) $f \in \dbp$. Part (ii) should be interpreted pointwise $\nu$-almost everywhere in $F$.
\end{remark}

\begin{definition}\label{de:partition}
{\rm (i)} Let $(\psi^Z_x)_{x \in X^Z}$ be a collection of Lipschitz functions $\psi^Z_x \colon Z \to [0,1]$ such that $\psi^Z_x$ is supported on $B(x)$ for all $x$, $\Lip \psi^Z_x \lesssim 2^{|x|}$ for all $|x|$ and $(\psi^Z_x)_{x \in X^Z_n}$ is a partition of unity of $Z$ for all $n \in \integer$. Define the collection of functions $(\psi^F_x)_{x \in X^F}$ in the same way with $F$ in place of $Z$.

\smallskip
{\rm (ii)} For $u\colon X^Z \to \complex$, define $T^Z_n u \colon Z \to \complex$ for all $n \in \integer$ by
\[
  T^Z_n u = \sum_{x \in X^Z_n} u(x) \psi^Z_x.
\]
Define $T^F_n u$ for $u \colon X^F\to \complex$ and $n \in \integer$ analogously.
\end{definition}

For $f \in \locint(Z)$ we obviously have $\lim_{n \to \infty}T^Z_n (Pf) = f$ pointwise $\mu$-almost everywhere (e.g.~at the Lebesgue points of $f$). It also turns out that for $f \in \dfp$ with suitable indices, $T^Z_n (Pf)$ approximates $f$ in the quasinorm of $\dfp$ as $n \to \infty$, and a similar result holds in the scale $\dbp$; see \cite[Theorem 3.3]{BSS} and \cite[Theorem 3.2]{S}. All this of course holds with $F$ in place of $Z$.

\begin{proof}[Proof of Theorem \ref{th:besov-trace}]
We will first construct the trace operator $\res$. In fact, as explained in Remark \ref{re:assumptions} above, we could take part (i) of the statement as the definition of $\res$, but we shall construct the operator in a slightly more roundabout way so that the boundedness becomes evident.

Letting $f \in \dbp$, we have $\|f\|_{\dbp} = \|d(Pf)\|_{\I^s_{p,q}(E^Z)}$, and consequently (by \eqref{eq:besov-norm-2} and the Ahlfors regularity of the spaces $Z$ and $F$),
\beqla{eq:derivative-norm}
  \big\| d(Pf)_{|E^F} \big\|_{\I^{s-\gamma/p}_{p,q}(E^F)} \approx \big\| d(Pf)_{|E^F} \big\|_{\I^{s}_{p,q}(E^Z)} \leq \|d(Pf)\|_{\I^s_{p,q}(E^Z)} = \|f\|_{\dbp} < \infty.
\eeq
We then need estimate $\| f_{|F} \|_{\db^{s-\gamma/p}_{p,q}(F)}$  in terms of the leftmost quantity above. To this end, write
\[
  I^F_n u := \sum_{(y,y')\in (X^F_n \times X^F_{n+1}),\, y\sim y'} u(e_{y,y'})\psi^F_y \psi^F_{y'}
\]
for all sequences $u$ defined on $E^F$ and integers $n \in \integer$, and fix $\xi_0 \in F$. We have that
\[
  \integral^F \big( d(Pf)_{|E^F} \big) := \lim_{N\to\infty} \bigg( \sum_{n=-N}^{N} I^F_n \big(d(Pf)_{|E^F}\big)(\cdot) - \sum_{n=-N}^{-1} I^F_n \big(d(Pf)_{|E^F}\big)(\xi_0) \bigg)
\]
converges in $\locint(F)$ and pointwise $\nu$-almost everywhere (see Lemma \ref{le:convergence} in the Appendix below). According to \cite[Proposition 4.3]{S} (see also \cite[Proposition 6.3]{BSS}), the $\db^{s-\gamma/p}_{p,q}(F)$-norm of the limit function is bounded from above by a constant times the leftmost quantity in \eqref{eq:derivative-norm}, and hence by a constant times $\|f\|_{\dbp}$.

Note that we have $I^F_n(du) = T^F_{n+1} u - T^F_n u$ for all sequences $u\colon X\to\complex$ (by definition), and since $d(Pf)_{|E^F}$ (as a sequence on $E^F$) is simply obtained as the discrete derivative of $(Pf)_{|X^F}$, we get
\beqla{eq:restriction-1}
  I^F_n\big( d(Pf)_{|E^F} \big) = T^F_{n+1} \big((Pf)_{|X^F}\big) - T^F_{n} \big((Pf)_{|X^F}\big)
\eeq
for all $n$. By Lemma \ref{le:convergence} in the Appendix below, we also have
\beqla{eq:restriction-2}
  \lim_{M \to -\infty} \Big( T^F_M\big((Pf)_{|X^F}\big)(\xi) - T^F_M\big((Pf)_{|X^F}\big)(\xi_0)\Big) = 0
\eeq
for all $\xi \in F$. Combining \eqref{eq:restriction-1} and \eqref{eq:restriction-2} with the fact that $\nu$-almost every point of $F$ is a Lebesgue point of $f$, we get
\begin{align*}
  \integral^F \big( d(Pf)_{|E^F} \big)(\xi) & = \lim_{N \to \infty} \Big( T^F_{N+1}\big( (Pf)_{|X^F}\big)(\xi) - T^F_{-N}\big( (Pf)_{|X^F}\big)(\xi) \\
& \qquad - T^F_{0}\big( (Pf)_{|X^F}\big)(\xi_0) + T^F_{-N}\big( (Pf)_{|X^F}\big)(\xi_0)  \Big)\\
  & = f_{|F}(\xi) - T^F_0\big((Pf)_{|X^F}\big)(\xi_0)
\end{align*}
for $\nu$-almost all $\xi \in F$, where $T^F_0\big((Pf)_{|X^F}\big)(\xi_0)$ is a constant. We can therefore take $\res f := \integral^F \big( d(Pf)_{|E^F} \big) + T^F_0\big((Pf)_{|X^F}\big)(\xi_0)$.

We shall next construct the extension operator $\ext$ with the additional assumption that $q < \infty$. In this case, it suffices to construct a bounded linear operator $\ext\colon \db^{s-\gamma/p}_{p,q}(F)\to \dbp$ satisfying (ii) for Lipschitz functions $f \in \dbp$ with bounded support, since these functions form a dense subspace of $\dbp$ \cite[Corollary 3.3]{S}. Taking $f \in \db^{s-\gamma/p}_{p,q}(F)$, $u := d(Pf)$ is a priori defined as a sequence on $E^F$, but it can be extended to $E^Z$ simply by defining $u(e) = 0$ for all $e \in E^Z\backslash E^F$. Now $\|u\|_{\I^s_{p,q}(E^Z)} \approx \|d(Pf)\|_{\I^{s-\gamma/p}_{p,q}(E^F)} < \infty$, so
\[
  \ext f := \integral^Z u + T^Z_0 (Pf)(\xi_0) = \lim_{N\to\infty} \bigg( \sum_{n=-N}^{N} I^Z_n u (\cdot) - \sum_{n=-N}^{-1} I^Z_n u(\xi_0) \bigg) + T^Z_0(Pf)(\xi_0),
\]
where $\xi_0$ is a fixed point of $F$, converges in $\locint(Z)$ and pointwise $\mu$-almost everywhere to a function in $\dbp$ with norm bounded from above by a constant times $\|f\|_{\db^{s-\gamma/p}_{p,q}(F)}$.

To verify the condition (ii) for Lipschitz functions $f \in \dbp$ with bounded support, we first show that in this case the series defining $\integral^Z u$ converges everywhere in $Z$ to a continuous function. By the Lipschitz continuity of $f$, we have $\sup_{\xi \in Z} |I^Z_n u(\xi)| \lesssim 2^{-n} \Lip (f)$ for all $n \in \integer$, so $\sum_{n \geq 0} I^Z_n u$ converges uniformly in $Z$. Furthermore, by the Lipschitz continuity of the functions $\psi^Z_x$, $x \in X^Z$, we have $|I_n^Z u(\xi) - I_n^Z u(\xi_0)| \lesssim 2^n d(\xi,\xi_0) \|f\|_{L^{\infty}(F)}$, so the series $\sum_{n < 0} \big(I_n^Z u(\cdot) - I_n^Z u(\xi_0)\big)$ converges uniformly on boun\-ded subsets of $Z$. All in all, the series defining $\integral^Z u$ converges uniformly on bounded subsets of $Z$, and the limit function must hence be continuous in $Z$.

Now since $u$ is not in general obtained as a discrete derivative of a sequence on $X^Z$, we do not have an analog of \eqref{eq:restriction-1} on $Z$. However, by the choices of the hyperbolic fillings (see Lemma \ref{le:filling-choice}), $B(x) \cap F = \emptyset$ for $x \in X^Z\backslash X^F$. Hence for all $\xi \in F$ and $n \in \integer$,
\begin{align*}
  I^Z_n u(\xi) & = \sum_{(y,y')\in (X^Z_n \times X^Z_{n+1}),\, y\sim y'} u(e_{y,y'})\psi^Z_y(\xi) \psi^Z_{y'}(\xi) \\
& = \sum_{(y,y')\in (X^F_n \times X^F_{n+1}),\, y\sim y'} u(e_{y,y'})\psi^Z_y(\xi) \psi^Z_{y'}(\xi) \\
& = \sum_{(y,y')\in (X^F_n \times X^F_{n+1}),\, y\sim y'} \big( Pf(y') - Pf(y)\big)\psi^Z_y(\xi) \psi^Z_{y'}(\xi) \\
& = T^Z_{n+1}(Pf)(\xi) - T^Z_{n}(Pf)(\xi).
\end{align*}
By continuity, all points of $Z$ are Lebesgue points of $\integral^Z u$ and all points of $F$ are Lebesgue points of $f$. Combining this with the formula above and Lemma \ref{le:convergence} in the Appendix below we get
\begin{align*}
  \integral^Z u (\xi) & = \lim_{N\to\infty}\Big( T^Z_{N+1}(Pf)(\xi) - T^Z_{-N}(Pf)(\xi) - T^Z_0(Pf)(\xi_0) + T^Z_{-N}(Pf)(\xi_0) \Big)  \\
& = f(\xi) - T^Z_0(Pf)(\xi_0)
\end{align*}
for $\xi \in F$; note that Lemma \ref{le:convergence} applies here since the Ahlfors regularity of $Z$ means that either $\diam(Z) < \infty$ or $\mu(Z) = \infty$. Altogether,
\[
  \res\big(\ext(f)\big) = \res\big(\integral^Z u\big) + \res\big( T^Z_0(Pf)(\xi_0) \big) = \big( f - T^Z_0(Pf)(\xi_0) \big) + T^Z_0(Pf)(\xi_0) = f
\]
$\nu$-almost everywhere in $F$.

Finally, if $q = \infty$, the operator $\ext$ constructed above extends to a bounded linear operator from $\db^{s-\gamma/p}_{p,\infty}(F)$ to $\db^s_{p,\infty}(Z)$ satisfying (ii), since $\db^{s-\gamma/p}_{p,\infty}(F)$ is obtained as a real interpolation space between the spaces $\db^{s_0-\gamma/p}_{p,p}(F)$ and $\db^{s_1-\gamma/p}_{p,p}(F)$ with $s_0 < s < s_1$ and $s_1 - s_0 \ll 1$, and similarly for the space $\db^s_{p,\infty}(Z)$ \cite[Theorem 4.3]{HIT}.
\end{proof}

\section{Traces of Triebel-Lizorkin and Sobolev spaces}\label{se:triebel}

In this section we shall give the proofs of Theorems \ref{th:triebel-trace} and \ref{th:sobolev-trace}. The assumptions on the metric measure spaces are as in the statements of these theorems -- $Z := (Z,d,\mu)$ is $Q$-Ahlfors regular, and $F$ is a closed and porous subset of $Z$ equipped with the metric $d_{|F}$ and the $\lambda$-Hausdorff measure $\nu$, which is assumed to be $\lambda$-Ahlfors regular. The hyperbolic fillings $(X^Z,E^Z)$ and $(X^F,E^F)$ are chosen as in the previous section.

The main observation concerning the porosity of $F$ is that now the $\J^s_{p,q}(E^Z)$-norm of a sequence living ``above'' $F$ is essentially independent of $q$. In \cite[Theorem 13.7]{FJ}, a similar phenomenon was observed for sequence spaces corresponding to the Fourier-analytically defined Triebel-Lizorkin spaces in the Euclidean setting, under a slightly weaker condition called \emph{NST} instead of porosity.

\lem{le:porosity}
Let $s \in (0,\infty)$, $p \in (0,\infty)$ and $q,\,q' \in (0,\infty]$. Then for all sequences $u\colon E^Z \to \complex$ supported on $E^F$ we have
\[
  \|u\|_{\J^s_{p,q}(E^Z)} \approx \|u\|_{\J^s_{p,q'}(E^Z)},
\]
with the implied constants independent of $u$.
\elem

\begin{proof}
The porosity of $F$ means that for all $e \in E^Z$ such that $B(e) \cap F \neq \emptyset$, we can take $x_e \in X^Z$ such that $|e| - \sigma \leq |x_e| \leq |e|$ for some fixed $\sigma \geq 0$ and $2B(x_e) \subset B(e) \backslash F$. By Remark \ref{re:space-properties} (vii) we thus have
\[
  \|u\|_{\J^s_{p,q}(E^Z)} \approx \bigg(\int_Z \bigg( \sum_{e \in E^Z\, : \, B(e)\cap F \neq \emptyset} \big[ 2^{|e|s} |u(e)| \big]^q \chi_{B(x_e)}(\xi) \bigg)^{p/q} d\mu(\xi) \bigg)^{1/p}
\]
(obvious modification for $q = \infty$), and now it suffices to show that for each $\xi \in Z$, only a uniformly finite number of terms in the latter sum are nonzero. Now suppose that $B(x_e)\cap B(x_{e'}) \neq \emptyset$ with $e$ and $e'$ like in the sum, and without loss of generality $|x_{e'}| \leq |x_e|$. By assumption we have $2^{\sigma'} B(x_e) \supset B(e)$ for some (universal) $\sigma'>0$,  so  $2^{\sigma'} B(x_e) \cap F \neq \emptyset$, which further means that $2^{\sigma'} B(x_e)$ is not contained in $2B(x_e')$. By construction this means that $|x_e| \leq |x_e'| + \sigma''$ for some uniform $\sigma'' \geq 0$. All in all, $\#\{e \in E^Z \,:\, B(e) \cap F \neq \emptyset \text{ and } B(x_e) \owns \xi \}$ is bounded uniformly in $\xi \in Z$, completing the proof.
\end{proof}

With this in mind we can give the proof of Theorem \ref{th:triebel-trace}. Remark \ref{re:assumptions}, with obvious modifications, holds here as well.

\begin{proof}[Proof of Theorem \ref{th:triebel-trace}]
Lemma \ref{le:porosity} tells us that for $f \in \dfp$, we have
\[
  \|f\|_{\dfp} \gtrsim \big\| d(Pf)_{|E^F} \big\|_{\J^s_{p,q}(E^Z)} \approx \big\| d(Pf)_{|E^F} \big\|_{\J^s_{p,p}(E^Z)} = \big\| d(Pf)_{|E^F} \big\|_{\I^s_{p,p}(E^Z)}.
\]
$\res$ can thus be constructed as in the proof of Theorem \ref{th:besov-trace}. To construct $\ext$, note that for all $f \in \db^{s-(Q-\lambda)/p}_{p,p}(F)$, the sequence $u := d(Pf)$ defined a priori on $E^F$ can be extended as zero on $E^Z \backslash E^F$, so by Lemma \ref{le:porosity} we have
\[
  \|f\|_{\db^{s-(Q-\lambda)/p}_{p,p}(F)} = \|d(Pf)\|_{\I^{s-(Q-\lambda)/p}_{p,p}(F)} \approx \|u\|_{\I^s_{p,p}(E^Z)} \approx \|u\|_{\J^s_{p,q}(E^Z)},
\]
and since Lipschitz functions with bounded support are dense in $\db^{s-(Q-\lambda)/p}_{p,p}(F)$, we may proceed as in the proof of Theorem \ref{th:besov-trace}.
\end{proof}

Theorem \ref{th:sobolev-trace} can be proven by utilizing the same idea. Again, Remark \ref{re:assumptions} with obvious modifications holds here as well.

\begin{proof}[Proof of Theorem \ref{th:sobolev-trace}]
We have $\dm^{1,p}(Z) = \df^1_{p,\infty}(Z)$ \cite[Proposition 3.1]{BSS}. The operator $\res$ can thus be constructed exactly as in the proof of Theorem \ref{th:triebel-trace} above -- here it is important remember that the porosity assumption on $F$ means that $\lambda < Q$ (see Remark \ref{re:porous}), which in turn implies $1 - (Q-\lambda)/p < 1$, which one needs to guarantee that the series defining $\integral^F\big(d(Pf)_{|E^F}\big)$ converges to a function in $\db^{1-(Q-\lambda)/p}_{p,p}(F)$ \cite[Proposition 4.3]{S}.

To construct the operator $\ext$, suppose that $f \in \db^{1-(Q-\lambda)/p}_{p,p}(F)$. Again, extending the sequence $u := d(Pf)$ as zero on $E^Z\backslash E^F$, we have
\[
  \|f\|_{\db^{1-\frac{Q-\lambda}{p}}_{p,p}(F)} \approx \|u\|_{\J^1_{p,p}(E^Z)} \approx \|u\|_{\J^1_{p,1}(E^Z)}.
\]
By Lemma \ref{le:convergence} in the Appendix below, we therefore have that
\[
  \integral^Z u := \lim_{N\to\infty} \bigg( \sum_{n=-N}^{N} I^Z_n u (\cdot) - \sum_{n=-N}^{-1} I^Z_n u(\xi_0) \bigg),
\]
where $\xi_0$ is a fixed point of $F$, converges in $\locint(Z)$ and pointwise $\mu$-almost everywhere (but \emph{not} to a function in $\df^{1}_{p,1}(Z)$, since \cite[Proposition 6.3]{BSS} requires $s < 1$). For all $\xi$ and $\eta$ outside of a set of zero $\mu$-measure we can however use the Lipschitz continuity of the functions $\psi^Z_y$ to obtain
\begin{align*}
  |\integral^Z u(\xi) - \integral^Z u(\eta)|
& \leq \sum_{n \in \integer} \sum_{(y,y')\in (X^Z_n \times X^Z_{n+1}),\, y\sim y'} |u(e_{y,y'})| \big| \psi^Z_y(\xi) \psi^Z_{y'}(\xi) - \psi^Z_y(\eta) \psi^Z_{y'}(\eta)\big| \\
& \lesssim d(\xi,\eta) \sum_{n \in \integer} 2^n \sum_{(y,y')\in (X^Z_n \times X^Z_{n+1}),\, y\sim y'} |u(e_{y,y'})| \big(\chi_{B(e_{y,y'})}(\xi) + \chi_{B(e_{y,y'})}(\eta)\big) \\
& \lesssim d(\xi,\eta)\bigg( \sum_{e \in E^Z} 2^{|e|} |u(e)| \chi_{B(e)}(\xi) + \sum_{e \in E^Z} 2^{|e|} |u(e)| \chi_{B(e)}(\eta) \bigg), 
\end{align*}
so
\[
  \|\integral^Z u \|_{\dm^{1,p}(Z)} \lesssim \bigg( \int_Z \Big( \sum_{e \in E^Z} 2^{|e|} |u(e)| \chi_{B(e)}(\xi) \Big)^p d\mu(\xi)\bigg)^{1/p} = \|u\|_{\J^1_{p,1}(E^Z)} \approx \|f\|_{\db^{1-\frac{Q-\lambda}{p}}_{p,p}(F)},
\]
and since Lipschitz functions with bounded support are dense in $\db^{1-(Q-\lambda)/p}_{p,p}(F)$, we may proceed as in the proof of Theorem \ref{th:besov-trace}.
\end{proof}

\section{Traces of non-homogeneous function spaces}\label{se:nonhomogeneous}

In this section we provide a definition of the non-homogeneous function spaces $\ndb^s_{p,q}$ and $\ndf^s_{p,q}$ in terms of the discrete derivatives and a fixed-level trace approximation of the Poisson extension $Pf$ of a function $f \in \locint$. This definition is of independent interest, and we shall apply it to prove Theorem \ref{th:nonhomogeneous}, i.e.~the non-homogeneous counterpart of our main theorems. A theory of non-homogeneous Besov and Triebel-Lizorkin spaces in the context of \emph{reverse doubling} metric measure spaces using spaces of test functions has previously been developed e.g.~in \cite{HMY}.

Our basic assumption here is that the space $(Z,d,\mu)$ satisfies the doubling condition \eqref{eq:doubling}. The operator $T^Z_0$ is as in Definition \ref{de:partition}.

\begin{definition}\label{de:nonhomogeneous-spaces}
\rm{(i)} Suppose that $0 < s \leq 1$, $Q/(Q+s) < p \leq \infty$ and $0 < q \leq \infty$. Then $\ndbp$ is the quasi-normed space of functions $f \in \locint(Z)$ such that
\[
  \|f\|_{\ndbp} := \| T^Z_0(Pf) \|_{L^p(Z)} + \Big(\sum_{k \geq 0}2^{ksq}\big\|\sum_{|e| = k} |d(Pf)(e)| \chi_{B(e)}\big\|_{L^p(Z)}^q\Big)^{1/q}
\]
(standard modification for $q = \infty$) is finite.

\smallskip
\rm{(ii)} Suppose that $0 < s \leq 1$, $Q/(Q+s) < p < \infty$ and $Q/(Q+s) < q \leq \infty$. Then $\ndfp$ is the quasi-normed space of functions $f \in \locint(Z)$ such that
\[
  \|f\|_{\ndfp} := \| T^Z_0(Pf) \|_{L^p(Z)} + \bigg( \int_Z \Big(\sum_{|e| \geq 0} \big[2^{|e|s}|d(Pf)(e)|\big]^q \chi_{B(e)}(\xi )\Big)^{p/q} d\mu(\xi) \bigg)^{1/p}
\]
(standard modification for $q = \infty$) is finite.
\end{definition}

Note that unlike in the homogeneous case, $\| \cdot \|_{\ndb^s_{p,q}}$ and $\|\cdot \|_{\ndf^s_{p,q}}$ are honest quasi-norms (modulo equality $\mu$-almost everywhere). Many properties of the homogeneous spaces also hold for the spaces $\ndb$ and $\ndf$, as evidenced by the alternative characterization of the non-homogeneous quasi-norms given in the Proposition below -- let us simply mention here that they are quasi-Banach spaces, and a standard argument shows that they are reflexive Banach spaces when $1 < p,\,q < \infty$. The Proposition below also shows that in the smoothness range $0 < s < 1$, these spaces coincide with the non-homogeneous spaces considered in e.g.~\cite{HIT} and \cite{HKT}, and with the standard Fourier-analytically defined non-homogeneous spaces when $Z = \real^d$. The Proposition further shows that we have $\ndf^1_{p,\infty}(Z) = M^{1,p}(Z)$ for $Q/(Q+1) < p < \infty$, where the space $M^{1,p}$ is as defined in \cite{H,H2}

\prop{pr:nonhomogeneous-norm}
\rm{(i)} Let $s$, $p$ and $q$ be admissible parameters for the space $\ndbp$. Then $\ndbp = L^p(Z)\cap \dbp$, and we have
\[
  \|f\|_{\ndbp} \approx \|f\|_{L^p(Z)} + \|f\|_{\dbp}
\]
for all $f \in \locint(Z)$, with the implied constants independent of $f$.

\smallskip
\rm{(ii)} Let $s$, $p$ and $q$ be admissible parameters for the space $\ndfp$. Then $\ndfp = L^p(Z)\cap \dfp$, and we have
\[
  \|f\|_{\ndfp} \approx \|f\|_{L^p(Z)} + \|f\|_{\dfp}
\]
for all $f \in \locint(Z)$, with the implied constants independent of $f$.
\eprop

Before proving this, let us formulate the following auxiliary result, which we shall use both in the proof of this Proposition as well as in the proof of Theorem \ref{th:nonhomogeneous} for the cases with $p < 1$. The proof of this result is given in the Appendix below. The Lemma essentially gives a local embedding of $\ndbp$ and $\ndfp$ into $L^1(Z)$ as long as $p > Q/(Q+s)$.

\lem{le:nonhomogeneous-small-p}
Suppose that $0 < p < 1$ and $\epsilon > 0$ is defined by $p = Q/(Q+\epsilon)$.

\smallskip
{\rm (i)} Let $x \in X^Z_0$. Then
\begin{align*}
  \big(\|f\|_{L^1(B(x))}\big)^p & \lesssim \mu\big( B(x) \big)^{p-1} \bigg( \int_{ B(x)} |T^Z_0(Pf)|^p d\mu \\
& \qquad + \sum_{k \geq 0} 2^{k\epsilon p} \big\|\sum_{|e| = k} |d(Pf)(e)| \chi_{B(e)}\big\|_{L^p(\sigma B(x))}^p \bigg)
\end{align*}
for all $f \in \locint(Z)$, where $\sigma \geq 1$ is a constant depending only on the choice of the hyperbolic filling $(X^Z,E^Z)$ (see Remark \ref{re:filling-choice}) and the implied constant is independent of $x$ and $f$.

\smallskip
{\rm (ii)} Let $B \subset Z$ be a ball with radius $r \geq 1$. Then
\begin{align*}
  \big(\|f\|_{L^1(B)}\big)^p & \lesssim r^{Q-Qp}\mu( B )^{p-1} \bigg( \int_{\sigma B} |T^Z_0(Pf)|^p d\mu \\
& \qquad + \sum_{k \geq 0} 2^{k\epsilon p} \big\|\sum_{|e| = k} |d(Pf)(e)| \chi_{B(e)}\big\|_{L^p( \sigma B)}^p \bigg)
\end{align*}
for all $f \in \locint(Z)$, where $\sigma \geq 1$ is a constant depending only on the choice of the hyperbolic filling, and the implied constant is independent of $r$, $B$ and $f$.
\elem

\begin{proof}[Proof of Proposition \ref{pr:nonhomogeneous-norm}]
We shall only consider the case of Besov spaces, since the case of Triebel-Lizorkin spaces can be handled in a similar manner.

Suppose that $f \in \locint(Z)$ and take $\varepsilon \in (0,s)$ (or $\varepsilon = 0$ if $p \leq 1$). Since $T^Z_k(Pf) \stackrel{k\to\infty}{\longrightarrow} f$ pointwise $\mu$-almost everywhere in $Z$, we get
\begin{align}\label{eq:nonhomogeneous-convergence}
  \|f - T^Z_0(Pf)\|_{L^p(Z)} & \lesssim \Big( \sum_{k \geq 0} 2^{k \varepsilon p}\| T^Z_{k+1}(Pf) - T^Z_k(Pf) \|_{L^p(Z)}^p\Big)^{1/p} \\
& \lesssim \Big( \sum_{k \geq 0} 2^{k \varepsilon p}\big\|\sum_{|e| = k} |d(Pf)(e)| \chi_{B(e)}\big\|_{L^p(Z)}^p \Big)^{1/p}\notag\\
& \lesssim \Big(\sum_{k \geq 0} 2^{k s q}\big\|\sum_{|e| = k} |d(Pf)(e)| \chi_{B(e)}\big\|_{L^p(Z)}^q \Big)^{1/q}\notag
\end{align}
(obvious modifications for $p = \infty$ and/or $q = \infty$). The latter quantity appears on both sides of the desired estimate, so it remains to verify that the part of $\|f\|_{\dbp}$ corresponding to edges at levels $k < 0$ is bounded by $\|f\|_{\ndbp}$. To this end, for each $e$ with $|e| < 0$ choose a ball $B_e$ that contains $B(e)$ and has radius uniformly comparable to $2^{-|e|}$ (without loss of generality we may also assume that this radius is always $\geq 1$). The quantity we wish to estimate is thus
\begin{align*}
  & \bigg( \sum_{k < 0} 2^{ksq} \big\| \sum_{|e| = k} |d(Pf)(e)| \chi_{B(e)} \big\|_{L^p(Z)}^q \bigg)^{1/q} \\
& \qquad \lesssim \bigg( \sum_{k < 0} 2^{ksq} \Big( \sum_{|e| = k} \mu(B_e)^{1-p}\|f\|_{L^1(B_e)}^p \Big)^{q/p} \bigg)^{1/q} =: I.
\end{align*}
If $p \geq 1$, the desired estimate follows easily by using H\"older's inequality for each term $\mu(B_e)^{1-p}\|f\|_{L^1(B_e)}^p$ and noting that the balls $B_e$ at each fixed level have uniformly bounded overlap:
\begin{align*}
I \lesssim \|f\|_{L^p(Z)} \leq \|f - T^Z_0(Pf)\|_{L^p(Z)} + \|T^Z_0(Pf)\|_{L^p(Z)} \lesssim \|f\|_{\ndbp}.
\end{align*}
If on the other hand $Q/(Q+s) < p < 1$, take $\epsilon \in (0,s)$ such that $p = Q/(Q+\epsilon)$. Then part (ii) of Lemma \ref{le:nonhomogeneous-small-p} above yields
\begin{align*}
  \sum_{|e| = k} \mu(B_e)^{1-p}\|f\|_{L^1(B_e)}^p & \lesssim 2^{-k(Q-Qp)} \Big( \int_Z |T^Z_0(Pf)|^p d\mu\\
& \qquad + \sum_{k \geq 0} 2^{k\epsilon p} \big\|\sum_{|e| = k} |d(Pf)(e)| \chi_{B(e)}\big\|_{L^p(Z)}^p\Big) \\
& = 2^{-k(Q-Qp)} \|f\|_{\ndb^\epsilon_{p,p}(Z)}^p
\end{align*}
for each $k < 0$. As $\epsilon < s$, it is easily seen that $\|f\|_{\ndb^\epsilon_{p,p}(Z)} \lesssim \|f\|_{\ndbp}$, and since $Q+s - Q/p > 0$, we arrive at
\[
  I \lesssim \Big( \sum_{k < 0} 2^{k(Q+s-Q/p)q}\Big)^{1/q} \|f\|_{\ndbp} \approx \|f\|_{\ndbp}.\qedhere
\]
\end{proof}

We will further need the following density result, which also is of independent interest. A similar result has been recently established in \cite{HKT}.

\prop{pr:nonhomogeneous-density}
Suppose that $0 < s < 1$, $Q/(Q+s) < p \leq \infty$ and $0 < q < \infty$.

\smallskip
{\rm (i)} If $f \in \ndbp$, then $T^Z_n(Pf) \to f$ in the quasinorm of $\ndbp$ as $n \to \infty$.

\smallskip
{\rm (ii)} Suppose that in addition $p < \infty$ and $q > Q/(Q+s)$. If $f \in \ndfp$, then $T^Z_n(Pf) \to f$ in the quasinorm of $\ndfp$ as $n \to \infty$.
\eprop

\begin{proof}
The analogous results for the homogeneous versions of these spaces can be found in \cite[Theorem 3.3]{BSS} and \cite[Theorem 3.2]{S}. By Proposition \ref{pr:nonhomogeneous-norm}, it thus suffices to verify that the sequence $(T^Z_n(Pf))_{n \geq 0}$ converges to $f$ in the quasinorm of $L^p(Z)$. For simplicity we shall only consider the situation of part (i). In fact, this is essentially contained in the estimate \eqref{eq:nonhomogeneous-convergence}: taking $\varepsilon \in (0,s)$, we get
\begin{align*}
  \|f - T^Z_n(Pf)\|_{L^p(Z)} & \lesssim \Big( \sum_{k \geq n} 2^{k \varepsilon p}\big\|\sum_{|e| = k} |d(Pf)(e)| \chi_{B(e)}\big\|_{L^p(Z)}^p \Big)^{1/p}\\
& \lesssim \Big(\sum_{k \geq n} 2^{k s q}\big\|\sum_{|e| = k} |d(Pf)(e)| \chi_{B(e)}\big\|_{L^p(Z)}^q \Big)^{1/q},
\end{align*}
and since $q < \infty$, the latter quantity tends to zero as $n \to \infty$.
\end{proof}

With Propositions \ref{pr:nonhomogeneous-norm} and \ref{pr:nonhomogeneous-density} in mind, we are in a position to give the proof the non-homogeneous versions of our main results.

\begin{proof}[Proof of Theorem \ref{th:nonhomogeneous}]
For brevity we shall only consider the traces and extensions of Besov functions, as the analogs for Theorems \ref{th:triebel-trace} and \ref{th:sobolev-trace} will then follow with more or less obvious modifications of the original proofs. The hyperbolic fillings $(X^Z,E^Z)$ and $(X^F,E^F)$ are again chosen as in the proof of Theorem \ref{th:besov-trace}, i.e~so that $X^F_n$ is in a natural way a subset of $X^Z_n$ for all $n \in \integer$. When necessary, we shall use the notations $P^Zf$, $P^Ff$, $B^Z(x)$, $B^F(x)$ etc.~with the obvious meanings to distinguish between the relevant operators and balls on different spaces. 

First, suppose that the parameters $s$, $p$ and $q$ are admissible and $f \in \ndbp$. We will construct the trace of $f$ to $F$ using a modified version of the operator $\res$ which is better suited for non-homogeneous case. Denote $E^F_+ := \{ e \in E^F \,:\, |e| \geq 0\}$ and write
\begin{align*}
  \nres f & := \integral^F \big( d(Pf)_{|E^F_+}\big) + T^F_0\big((Pf)_{|X^F}\big) = \sum_{n=0}^{\infty} I^F_n \big(d(Pf)_{|E^F}\big) + T^F_0\big((Pf)_{|X^F}\big) \\
& = \sum_{n=0}^{\infty} \Big( T^F_{n+1}\big((Pf)_{|X^F}\big) - T^F_n\big((Pf)_{|X^F}\big) \Big) + T^F_0\big((Pf)_{|X^F}\big).
\end{align*}
As before, the series defining $\integral^F \big( d(Pf)_{|E^F_+}\big)$ converges $\nu$-almost everywhere in $F$, and by Proposition \ref{pr:nonhomogeneous-norm} the limit function has $\db^{s-(Q-\lambda)/p}_{p,q}(F)$-quasinorm bounded by a constant times
\beqla{eq:derivative-norm-2}
  \big\| d(Pf)_{|E^F_+} \big\|_{\I^{s-(Q-\lambda)/p}_{p,q}(E^F)} \lesssim \|d(Pf)\|_{\I^s_{p,q}(E^Z)} = \|f\|_{\dbp} \lesssim \|f\|_{\ndbp}.
\eeq
We further have $\integral^F \big( d(Pf)_{|E^F_+}\big) = f_{|F} - T^F_0\big( (Pf)_{|X^F}\big)$, and hence $\nres f = f_{|F}$, $\nu$-almost everywhere in $F$. We thus have to show that
\beqla{eq:nonhomogeneous-trace-1}
  \big\|\integral^F \big( d(Pf)_{|E^F_+}\big)\big\|_{L^p(F)} \lesssim \| f \|_{\ndbp}
\eeq
and that
\beqla{eq:nonhomogeneous-trace-2}
  \big\| T^F_0\big( (Pf)_{|X^F}\big) \big\|_{L^P(F)} + \big\| T^F_0\big( (Pf)_{|X^F}\big) \big\|_{\db^{s-(Q-\lambda)/p}_{p,q}(F)} \lesssim \| f \|_{\ndbp}.
\eeq

For \eqref{eq:nonhomogeneous-trace-1}, take $\varepsilon \in (0,s - (Q-\lambda)/p)$ (or $\varepsilon = 0$ if $p \leq 1$). Then for $\nu$-almost all $\xi \in F$ we have
\begin{align*}
  \big| \integral^F \big( d(Pf)_{|E^F_+}\big)(\xi) \big| & \leq \sum_{n \geq 0} \sum_{(y,y')\in (X^F_n \times X^F_{n+1}),\, y\sim y'} |P^Zf(y) - P^Zf(y')|\chi_{B^F(y)}(\xi) \chi_{B^F(y')}(\xi) \\
& \lesssim \Big( \sum_{n \geq 0} 2^{n\varepsilon p} \sum_{e \in E^F, \, |e| = n} |d(P^Zf)(e)|^p \chi_{B^F(e)}(\xi) \Big)^{1/p}.
\end{align*}
The Ahlfors regularity of $Z$ and $F$ implies that
\beqla{eq:ahlfors-codimension-2}
  \nu\big( B^F(e)\big) \approx 2^{|e|(Q - \lambda)} \mu\big( B^Z(e)\big)
\eeq
for all $e \in E^F_+$, so we further get
\begin{align*}
  \big\| \integral^F \big( d(Pf)_{|E^F_+}\big) \big\|_{L^p(F)} & \lesssim \Big( \sum_{n \geq 0} 2^{n\varepsilon p} \sum_{e \in E^F, \, |e| = n} \nu\big( B^F(e)\big) |d(P^Zf)(e)|^p  \Big)^{1/p} \\
& \lesssim \Big( \sum_{n \geq 0} 2^{n(\varepsilon p+ Q - \lambda)} \sum_{e \in E^F, \, |e| = n} \mu\big( B^Z(e)\big) |d(P^Zf)(e)|^p  \Big)^{1/p},
\end{align*}
and since $\varepsilon p+ Q - \lambda < s p$, the latter quantity is dominated by the right-hand side of \eqref{eq:nonhomogeneous-trace-1}.

We next estimate the $L^p(F)$-quasinorm on the left-hand side of \eqref{eq:nonhomogeneous-trace-2}. Again since $\nu(B^F(x))$ is comparable to $\mu(B^Z(x))$ for all $x \in X^F_0$, we get
\begin{align}
  \big\| T^F_0\big( (Pf)_{|X^F}\big) \big\|_{L^P(F)} & \lesssim \Big(\sum_{x \in X^F_0} \nu\big(B^F(x)\big) |P^Zf(x)|^p \Big)^{1/p} \label{eq:nonhomogeneous-trace-3} \\
& \lesssim \Big(\sum_{x \in X^F_0} \mu\big(B^Z(x)\big)^{1-p} \|f\|_{L^1(B^Z(x),\mu)}^p \Big)^{1/p}.\notag
\end{align}
If $p \geq 1$, the latter quantity is easily estimated by $\|f\|_{L^p(Z)}$, and if $Q/(Q+s) < p < 1$, we can use part (i) of Lemma \ref{le:nonhomogeneous-small-p} to get an estimate in terms of $\|f\|_{\ndbp}$.

Finally to estimate the $\db^{s-(Q-\lambda)/p}_{p,q}(F)$-quasinorm in \eqref{eq:nonhomogeneous-trace-2}, note that the Lipschitz continuity of the functions $(\psi^F_x)_{x \in X^F_0}$ implies that
\[
  \big| T^F_0\big( (Pf)_{|X^F}\big)(\xi) - T^F_0\big( (Pf)_{|X^F}\big)(\eta) \big| \lesssim \min\big(1,d(\xi,\eta)\big)\sum_{x \in X^F_0} |P^Zf(x)|\big( \chi_{B^F_x}(\xi) + \chi_{B^F_x}(\eta) \big),
\]
so a fairly straightforward computation -- which we leave to the reader\footnote{A similar computation also works for the spaces $\dfp$ with either $s < 1$ or $s = 1$ and $q = \infty$.} -- gives an estimate of the $\db^{s-(Q-\lambda)/p}_{p,q}$-quasinorm of $T^F_0(Pf_{|X^F})$ in terms of the $L^p$-quasinorm of the function
\[
  \sum_{x \in X^F_0} |P^Zf(x)|\chi_{B^F_x},
\]
which is essentially the quantity on the right-hand side of \eqref{eq:nonhomogeneous-trace-3}.

Now to construct the extension of a function $f \in \ndb^{s-(Q-\lambda)/p}_{p,q}(F)$, we extend the sequences $u := d(Pf)$ and $v := Pf$ as zero on $E^Z\setminus E^F$ and $X^Z\setminus X^F$ respectively and put
\[
 \next f := \integral^Z \big( u_{|E^Z_+}\big) + T^F_0(v) = \sum_{n=0}^{\infty} I^Z_n \big(u_{|E^Z}\big) + T^F_0(v).
\]
That $\next f$ converges in $\locint(Z)$ and the operator $\next$ is a bounded from $\ndb^{s-(Q-\lambda)/p}_{p,q}(F)$ into $\ndbp$ follows easily from an appropriate reformulation of the estimates established  above for the operator $\nres$.

To verify that $\nres ( \next f ) = f$ for $f \in \ndb^{s-(Q-\lambda)/p}_{p,q}(F)$, we first consider the case with $q < \infty$. Then by Proposition \ref{pr:nonhomogeneous-density}, it suffices to consider functions $f$ that are Lipschitz continuous on every bounded subset of $F$. In this case it is easily verified that the series defining $\next f$ converges uniformly on bounded subsets of $Z$, so the limit function is continuous on $Z$, and hence we plainly have $\next f(\xi) = f(\xi)$ for all $\xi \in F$. For the case $q = \infty$, one concludes again by using interpolation \cite[Theorem 4.3]{HIT}.
\end{proof}

\section{Appendix}\label{se:appendix}

In this section we shall make some additional remarks concerning the main results of this paper and elaborate on some less interesting technicalities that were used in their proofs.

\subsection{On the regularity and dimensions of $Z$ and $F$}

In Theorems \ref{th:besov-trace} through \ref{th:nonhomogeneous}, we assumed that the metric measure space $(Z,d,\mu)$ and $(F,d_{|F},\nu)$ were respectively $Q$-Ahlfors regular and $\lambda$-Ahlfors regular. Upon examining the proofs, it is easy to see that the only essential way this assumption is used is in the estimates \eqref{eq:derivative-norm} and \eqref{eq:ahlfors-codimension-2}, i.e.~to guarantee that
\beqla{eq:ahlfors-codimension-1}
  \frac{\mu\big(B_Z(\xi,r)\big)}{\nu\big(B_F(\xi,r)\big)} \approx r^{\gamma}
\eeq
for all $\xi \in F$ and $0 < r < \diam(F)$, where $\gamma = Q - \lambda$. Therefore this kind of a condition, known in previous literature (see e.g.~\cite{MSS}) as \emph{Ahlfors co-dimension $\gamma$ regularity} of $F$ with respect to $Z$, suffices for our main results.

To elaborate on this, suppose that the space $(Z,d,\mu)$ satisfies the doubling condition \eqref{eq:doubling}, and $\nu$ is merely a doubling measure on the metric space $(Z,d_{|F})$ satisfying \eqref{eq:ahlfors-codimension-1} for some $\gamma \in [0,Q)$. An easy calculation then shows that
\[
  \nu\big(B_F(\xi,\lambda r) \big) \lesssim \lambda^{Q-\gamma} \nu \big(B_F(\xi,r) \big)
\]
for all $\xi \in F$, $r > 0$ and $\lambda \geq 1$, i.e.~$\lambda := Q - \gamma$ is an upper bound for the dimension of the space $(F,d_{|F},\nu)$. We make the further assumptions that $\mu(Z) = \infty$ if $Z$ is unbounded, and $\gamma > 0$ for Theorems \ref{th:triebel-trace}, \ref{th:sobolev-trace} and their non-homogeneous counterparts.

With these assumptions and notations, Theorems \ref{th:besov-trace}, \ref{th:triebel-trace}, \ref{th:sobolev-trace} and \ref{th:nonhomogeneous} continue to hold. However, Remark \ref{re:assumptions} does not hold, since we have no precise information on the Hausdorff dimension of $F$. Therefore part (i) of each Theorem has to be interpreted as originally stated.

\subsection{Convergence of the operator $\integral$}

Here we shall establish the $\locint$ convergence of the operator $\integral$, as well as the formula \eqref{eq:restriction-2}. The assumption for the Lemma below is that the space $(Z,d,\mu)$ satisfies the doubling property \eqref{eq:doubling}, and the notation is as in Definition \ref{de:partition}.

\lem{le:convergence}
{\rm (i)} Suppose that $u \in \big(\ipe\big) \cup \big(\jpe\big)$ with $0 < s < 1$, $Q/(Q+s) < p < \infty$ and $0 < q \leq \infty$, or $u \in \J^1_{p,q}(E)$ with $Q/(Q+1) < p < \infty$ and $0 < q \leq 1$. Then the limit
\[
  \integral^Z u : = \lim_{N\to\infty} \bigg( \sum_{n=-N}^{N} I^Z_n u (\cdot) - \sum_{n=-N}^{-1} I^Z_n u(\xi_0) \bigg),
\]
where $\xi_0 \in Z$ is a fixed point, exists in $\locint(Z)$ and pointwise $\mu$-almost everywhere. 

\smallskip
{\rm (ii)} Suppose that $u\colon X^Z \to \complex$ is a sequence such that $du \in \big(\ipe\big) \cup \big(\jpe\big)$ with $0 < s < 1$, $0 < p < \infty$ and $0 < q \leq \infty$. Then
\[
  \lim_{M \to -\infty} \Big( T^Z_Mu(\xi) - T^Z_Mu(\eta)\Big) = 0
\]
for all $\xi$, $\eta \in Z$. If either $\diam(Z) < \infty$ or $\mu(Z) = \infty$, the same holds for $s = 1$.
\elem

\begin{proof}
{\rm (i)} This is essentially included in the proofs of \cite[Proposition 6.3]{BSS} and \cite[Proposition 4.3]{S}, but since those results are only formulated in the case $0 < s < 1$, we shall repeat the main points of the argument here.

First, to obtain the convergence of the series $\sum_{n \geq 0} I^Z_n u$, fix $\epsilon \in (0,s)$ such that $r := Q/(Q+\epsilon) < p$. For a fixed ball $B \subset Z$ with radius $1$, we may use the doubling property \eqref{eq:doubling} as in the proof of \cite[Lemma 2.3]{BSS} to obtain
\[
  \sum_{n \geq 0} \int_B |I^Z_n u| d\mu \lesssim \sum_{\substack{|e| \geq 0 \\ B(e) \cap B \neq \emptyset}} \mu\big(B(e)\big) |u(e)| \lesssim \bigg(\sum_{\substack{|e| \geq 0 \\ B(e) \cap B \neq \emptyset}} \mu\big( B(e )\big)\big[2^{|e|\epsilon} |u(e)|\big]^r \bigg)^{1/r},
\]
where the implied constant in the last inequality depends on $B$. The latter quantity is essentially a local $I^\epsilon_{r,r}$ norm of $u$. Since we are only dealing with vertices $e$ with $|e| \geq 0$ and $\epsilon < s$, we may easily estimate this quantity from above by a local $\I^s_{r,q}$ or $\J^s_{r,q}$ norm (depending on which space $u$ belongs to), and since $r < p$, we obtain a finite upper bound by using H\"older's inequality. As the ball $B$ with radius $1$ was arbitrary, we have the desired convergence of $\sum_{n \geq 0}I^Z_n u$.

For the part corresponding to negative indices, we shall consider the case $u \in \jpe$ with admissible parameters; the case with $u \in \ipe$ can be handled in a similar manner. We shall examine the convergence of this part in $B := B(\xi_0,2^{-k})$, where $k < 0$ is an arbitrary integer. For all $n \leq k$, the Lipschitz continuity of the functions $\psi^Z_x$ yields
\[
  \int_{B}\big| I^Z_n u(\xi) - I^Z_n u(\xi_0)\big| d\mu(\xi) \lesssim \mu(B)2^{n-k} \sum_{\substack{|e| = n \\ B(e) \cap B \neq \emptyset}}|u(e)|,
\]
so that
\[
  \sum_{n \leq k} \int_{B}\big| I^Z_n u(\xi) - I^Z_n u(\xi_0)\big| d\mu(\xi) \lesssim \sum_{\substack{|e| \leq k\\B(e)\cap B \neq \emptyset}} 2^{|e|} |u(e)| \lesssim \Big( \sum_{\substack{|e| \leq k\\B(e)\cap B \neq \emptyset}} \big[ 2^{|e|s} |u(e)|\big]^q \Big)^{1/q}
\]
(obvious modification for $q = \infty$); here we used the knowledge that $\#\{ e \in E^Z \,:\, |e| = n \text{ and } B(e) \cap B \neq \emptyset\}$ is bounded uniformly in $n \leq k$ and the assumption that either $s < 1$ or $q \leq 1$. Writing $\lambda B(e) := \lambda B(e_-) \cup \lambda B(e_+)$ for a suitable $\lambda > 1$, we further get
\[
  \sum_{n \leq k} \int_{B}\big| I^Z_n u(\xi) - I^Z_n u(\xi_0)\big| d\mu(\xi) \lesssim \inf_{\xi \in B} \Big( \sum_{|e| \leq k} \big[ 2^{|e|s} |u(e)|\big]^q \chi_{\lambda B(e)}(\xi) \Big)^{1/q}.
\]
The latter infimum is finite, since by Remark \ref{re:space-properties} (vii), the function in question belongs to $L^p(Z)$.

\smallskip
(ii) There is nothing to prove if $\diam (Z) < \infty$, because then $T^Z_{M}u$ is simply a constant function for $M < 0$ with large enough absolute value. We can thus assume that $\diam(Z) = \infty$ in the remainder of this proof.

Write $\|du\|$ for $\|du\|_{\ipe}$ or $\|du\|_{\jpe}$, whichever is finite. With $\xi$ and $\eta$ fixed, we can for all $M < 0$ with sufficiently large absolute value find a vertex $x_M \in X^Z$ such that $B(x)$ contains both $\xi$ and $\eta$. By the doubling property, $\mu\big(B(x_M)\big)$ is comparable to $\mu\big(B(\eta,2^{-M})\big)$ uniformly in $M$. Using the Lipschitz continuity of the functions $\psi^Z_x$, we obtain
\begin{align*}
  \big| T^Z_Mu(\xi) - T^Z_Mu(\eta)\big| & \leq \sum_{x \in X^Z_{M}}\big| u(x) - u(x_M)|\big|\psi^Z_x(\xi) - \psi^Z_x(\eta) \big| \lesssim 2^{M} \sum_{e \,:\, x_M \in \{e_-,e_+\} }|(du)(e)| \\
  & \lesssim 2^{M(1-s)} \mu\big( B(x_M) \big)^{-1/p} \|du\| \approx 2^{M(1-s)} \mu\big( B(\eta,2^{-M}) \big)^{-1/p} \|du\|.
\end{align*}
Since either $s < 1$ or $\mu(Z) = \infty$, the latter quantity tends to zero as $M \to -\infty$.
\end{proof}

\subsection{Proof of Lemma \ref{le:nonhomogeneous-small-p}}

Finally we present the proof of the auxiliary result that was used in Section \ref{se:nonhomogeneous} when handling the cases with $p < 1$.

\begin{proof}[Proof of Lemma \ref{le:nonhomogeneous-small-p}]
(i) Note that
\beqla{eq:hv-estimate-1}
  \big(\|f\|_{L^1(B(x))}\big)^p \leq \|f - Pf(x)\|_{L^1(B(x))}^p + \mu\big( B(x)\big)^p |Pf(x)|^p.
\eeq
For the latter term, we have the following simple estimates:
\begin{align*}
  |Pf(x)|^p &= \dashint_{B(x)} |Pf(x)|^p d\mu(\xi) \\
& \leq \dashint_{B(x)} |Pf(x) - T^Z_0(Pf)(\xi)|^p d\mu(\xi) + \dashint_{B(x)} |T^Z_0(Pf)(\xi)|^p d\mu(\xi) \\
& \lesssim \mu\big( B(x) \big)^{-1} \Big( \big\|\sum_{|e| = 0} |d(Pf)(e)| \chi_{B(e)}\big\|_{L^p(\sigma B(x))}^p + \int_{B(x)} |T^Z_0(Pf)(\xi)|^p d\mu(\xi) \Big),
\end{align*}
and multiplying both sides by $\mu( B(x) )^p$ gives a good enough estimate for the latter term in \eqref{eq:hv-estimate-1}.

The necessary estimate for the term $\|f - Pf(x)\|_{L^1(B(x))}^p$ in \eqref{eq:hv-estimate-1} is basically contained in the proof of Lemma \ref{le:convergence} above, so we only give a rough outline here that takes into account the multiplicative constants depending on $\mu(B(x))$. Using again the fact that $T^Z_k(Pf) \to f$ in $\locint(Z)$ as $k \to \infty$, the assumption that $p = 1 - (\epsilon/Q)p$ and the doubling property of $\mu$, we get
\begin{align*}
  \|f - Pf(x)\|_{L^1(B(x))}^p & \lesssim \sum_{\substack{|e| \geq 0 \\ B(e)\cap B(x)\neq \emptyset}} \mu\big( B(e) \big)^p |d(Pf)(e)|^p \\
 & \lesssim \mu\big(B(x)\big)^{-(\epsilon/Q)p} \sum_{\substack{|e| \geq 0 \\ B(e)\cap B(x)\neq \emptyset}} 2^{|e|\epsilon p}\mu\big( B(e) \big) |d(Pf)(e)|^p \\
 & \lesssim \mu\big( B(x) \big)^{p-1} \sum_{k\geq 0} 2^{k\epsilon p}\big\|\sum_{|e| = k} |d(Pf)(e)| \chi_{B(e)}\big\|_{L^p(\sigma B(x))}^p,
\end{align*}
which finishes the proof of part (i).

\smallskip
(ii) This follows from part (i) by covering $B$ with an optimal collection of balls $B(x_i)$, $x_i \in X^Z_0$, and noting that the doubling property yields $\mu( B(x_i) )^{p-1} \lesssim r^{Q-Qp} \mu(B)^{p-1}$ for all $x_i$.
\end{proof}

\end{document}